\documentstyle[12pt,amssymbols]{article}
\newcommand{\bq}{\begin{equation}} 
\newcommand{\eq}{\end{equation}} 
\newcommand{\bqa}{\begin{eqnarray}} 
\newcommand{\eqa}{\end{eqnarray}} 
\newcommand{\ra}{\rightarrow}

\def\e{\eta}

\def\ed{\end{document}} 
 
\def\ra{\rightarrow} 
 
\def\2pi{1\over 2\pi i} 

\def\~{\tilde} 
\def\newline{\hfil\break}

\def\ra{\rightarrow}

\def\sq2{\sqrt{2}} 
\def\sqk2{\sqrt{2(k+2}} 
\def\sqk{\sqrt{k}}

\def\be{\begin{equation}} 
\def\ee{\end{equation}} 
\def\br{\begin{array}} 
\def\er{\end{array}} 
\def\bea{\begin{eqnarray}} 
\def\eea{\end{eqnarray}} 
\def\ba{\begin{equation}\begin{array}} 
\def\ea{\end{array}\end{equation}} 
\def\bac{\begin{equation}\begin{array}{rll}}

\newcommand{\ugl}{U_q (\widehat{gl(2)})}

\def\Z{{\Bbb Z}} 
\def\C{{\Bbb C}}

\begin{document} 
\rightline{ITP-SB-97-01}  
\rightline{January, 1997}  
\vbox{\vspace{-10mm}}  
\vspace{1.0truecm}  
\begin{center}  
{\LARGE \bf   
A Hopf algebra isomorphism between 
two realizations 
of the quantum affine algebra $U_q(\widehat{gl(2)})$ }\\[8mm]  
{\large A.H. Bougourzi$^{1}$ and  A. Sebbar$^{2}$  
}\\  
[6mm]{\it 
$^{1}$
Institute of Theoretical Physics\\  
SUNY at Stony Brook\\ 
Stony Brook, NY 11794
\\[3mm] 
$^{2}$ Department of Mathematics\\
SUNY at Stony Brook\\
Stony Brook, NY 11794
}\\[20mm]  
\end{center}  
\vspace{1.0truecm}  
\begin{abstract}
We consider the algebra isomorphism found by Frenkel and 
Ding between the $RLL$
and the Drinfeld realizations of 
$U_q(\widehat{gl(2)})$. After we note that  this is not a Hopf
algebra isomorphism, we prove that there is a unique
Hopf algebra 
structure for the Drinfeld realization so that this 
isomorphism becomes a Hopf
algebra isomorphism. Though more complicated, 
this Hopf algebra structure
is also closed, just as the one found previously by Drinfeld.

\end{abstract} 
\newpage 

\section{Introduction}  

Besides their rich mathematical structures, 
quantum affine algebras play a crucial role in the exact
calculation of some experimental 
physical quantities that have remained so far
unaccessible through traditional methods \cite{Boual96,Abaal96}. 
However, they are defined through various approaches and
therefore it is important to find the relations among all
of them.
In particular, 
it is well known that quantum affine algebras admit two
apparently different realizations in terms of currents. 
The first one which is based 
on the $RLL$ formalism was derived by Reshetikhin and 
Semenov-Tian-Shansky \cite{ReSe90}, 
following the general Faddeev-Reshtikhin-
Takhtajan matrix construction of quantum groups \cite{Fadal89}. 
The second one 
was derived by Drinfeld as a  quantum version  of  
 central extensions
of loop affine algebras \cite{Dri86}. 
The natural question  of 
finding the exact
connection between the two arises then. 
Frenkel and Ding addressed this question and were able to
construct through
the  Gauss  decomposition an algebra isomorphism 
between them \cite{DiFr93}. 
The crucial question 
of finding an isomorphism at
the level of Hopf algebras remained unanswered.
In fact, as we will see,  for the Hopf 
algebras structures given in \cite{DiFr93}, 
the above isomorphism is not 
a Hopf algebra isomorphism. 
However, we show that this algebra isomorphism can be 
upgraded to a Hopf algebra isomorphism 
for a different but unique Hopf algebra structure of the
Drinfeld realization. This  provides a second
example of a closed Hopf algebra structure for the latter
realization.
  
\section{Brief review of the Frenkel-Ding isomorphism}

Here we briefly recall the two realizations of $\ugl$ and the
algebra isomorphism between them as found by Frenkel and Ding 
\cite{DiFr93}. 
\subsection{The Reshetikhin-Semenov-Tian-Shansky realization
of $\ugl$}

Since this realization is based on a matrix representation of
quantum groups, and in particular on the $R(z)$ matrix which 
satisfies
the Yang-Baxter equation with a spectral parameter $z$, 
let us 
denote this algebra by $U(R(z))$. We review this realization 
following Ref. \cite{DiFr93} where a more compact but 
equivalent definition is
given.

{\bf Definition 2.1} \cite{ReSe90,DiFr93} $U(R(z))$ is an 
associative
algebra with a unit and is generated by the modes of the 
currents 
$\ell^\pm_{ij}(z)=\sum_{n=0}^\infty \ell^\pm_{ij}[\mp
n]z^{\pm n}, 
i,j=1,2$; with defining relations:
\bac
&&\ell^+_{21}[0]=\ell^-_{12}[0]=0\\
&&\ell^+_{ii}[0]\ell^-_{ii}[0]=\ell^-_{ii}[0]\ell^+_{ii}[0]=1,
\quad i=1,2;\\
&&R({z\over w})L^{\pm}_1(z)L^{\pm}_2(w)= 
L_2^{\pm}(w)L_1^{\pm}(z)R({z\over w}), \\ 
&&R({z_{+}\over w_{-}})L^{+}_1(z)L^{-}_2(w)= 
L_2^{-}(w)L_1^{+}(z)R({z_{-}\over w_{+}}). 
\label{Rel1}
\ea 
Here the notations are as follows:
 $L^\pm(z)=(\ell^\pm_{ij}(z))_{i,j=1}^2$, 
$L_1^\pm(z)=L^\pm(z)\otimes 1$, 
$L_2^\pm(z)=1\otimes L^\pm(z)$,  
$z_{\pm}=zq^{\pm c/2}$ with $q^c$ being the central element, and
\bac 
R(z)= 
{\scriptsize\left(\begin{array}{llll}1&0&0&0\\ 
0&{(1-z)q\over 1-zq^2}&{z(1-q^2)\over 
1-zq^2}&0\\
0&{1-q^2\over 
1-zq^2}&{(1-z)q\over 1-zq^2}&0\\ 
0&0&0&1\end{array} 
\right)}. 
\label{Rel2}
\ea 
Note also that in the above defining relations, $R(z)$ has a
Taylor expansion in $z$ (not in $z^{-1}$).  The important 
feature 
of this algebra is that it is equipped with a Hopf algebra
structure given by
\bac 
&&\Delta(L^{\pm}(z))=L^{\pm}(zq^{\pm (1\otimes
c/2)})\otimes^. 
L^{\pm}(zq^{\mp (c/2\otimes 1)}),\\
&&S(L^{\pm}(z))=L^{\pm}(z)^{-1},\\
&&{\cal E}(L^{\pm}(z))=I,
\label{Rel3}
\ea
where $\Delta$, $S$ and ${\cal E}$ denote the comultiplication,
the antipode and the counit respectively.

\subsection{\bf The Drinfeld realization of $\ugl$}

Here we briefly recall the defining relations of 
the Hopf algebra $\ugl$ as a quantum version of the central
extension of the loop algebra $\widehat{gl(2)}$.

{\bf Definition 2.2}\cite{Dri86,DiFr93}.
This is an associative algebra with a unit  and is generated by
the modes of the currents 
$X^{\pm}(z)=\sum_{n\in \Z}X^\pm_nz^{-n}$ and
 $k^\pm_i(z)=\sum_{n=0}^\infty k^\pm_{i,\mp n}z^{\pm n}$, 
$i=1,2$; and 
the central elements $q^{\pm {1\over 2}c}$, 
with defining relations:
\bac
&&k^+_{i0}k^-_{i0}=k^-_{i0}k^+_{i0}=1,\\
&&k^\pm_{i}(z)k^\pm_{i}(w)=k^\pm_{i}(w)k^\pm_{i}(z),\\
&&k^+_{i}(z)k^-_{i}(w)=k^-_{i}(w)k^+_{i}(z),\\  
&&{z_{\mp}-w_{\pm}\over z_{\mp} q-w_{\pm}q^{-1}}
k^\mp_{1}(z)k^\pm_{2}(w)=k^\pm_{2}(w)k^\mp_{1}(z)
{z_{\pm}-w_{\mp}\over z_{\pm} q-w_{\mp}q^{-1}},\\
&&k^\pm_{1}(z)X^+(w)k^\pm_{1}(z)^{-1}=
{z_{\pm}q-w q^{-1}\over z_{\pm} -w}X^+(w),\\
&&k^\pm_{2}(z)X^+(w)k^\pm_{2}(z)^{-1}=
{z_{\pm}q^{-1}-w q\over z_{\pm} -w}X^+(w),\\
&&k^\pm_{1}(z)^{-1}X^-(w)k^\pm_{1}(z)=
{z_{\mp}q-w q^{-1}\over z_{\mp} -w}X^-(w),\\
&&k^\pm_{2}(z)^{-1}X^-(w)k^\pm_{2}(z)=
{z_{\mp}q^{-1}-w q\over z_{\mp} -w}X^-(w),\\
&&(zq^{\pm 1}-wq^{\mp 1})X^\pm(z)X^\pm(w)=)X^\pm(w)X^\pm(z)
(zq^{\mp 1}-wq^{\pm 1}),\\
&&{[X^+(z),X^-(w)]}=(q-q^{-1})(\delta(z/wq^c)k^-_2(w_+)
k^-_1(w_+)^{-1}-\delta(zq^c/w)k^+_2(z_+)
k^+_1(z_+)^{-1}),
\label{Rel4}
\ea
where the formal $\delta$ function is defined by $\delta(z)=
\sum_{n\in \Z}z^n$.
This algebra is
equipped with a Hopf algebra structure. 
The
comultiplication, the antipode and the counit are given by 
the following relations:
\bac
\Delta(k^+_i(z))&=&k^+_i(zq^{-1\otimes {c\over 2}})\otimes
k^+_i(zq^{{c\over 2}\otimes 1}),\\
\Delta(k^-_i(z))&=&k^-_i(zq^{1\otimes {c\over 2}})\otimes
k^-_i(zq^{-{c\over 2}\otimes 1}),\\
\Delta(X^+(z))&=&X^+(z)\otimes 1+ 
k^-_{2}(zq^{c\otimes 1})
k^-_{1}(zq^{c\otimes 1})^{-1}
\otimes X^+(zq^{{c\over 2}\otimes 1}),\\
\Delta(X^-(z))&=&1\otimes X^-(z)+ 
 X^-(zq^{1\otimes c })\otimes 
k^+_{2}(zq^{1\otimes {c\over 2}})
k^+_{1}(zq^{1\otimes {c\over 2}})^{-1},\\
S(k^\pm_i(z))&=&k^\pm_i(z)^{-1},\\
S(X^+(z))&=&-k^-_{1}(z)k^-_{2}(z)^{-1}X^+(z),\\
S(X^-(z))&=&-X^-(z) k^+_{1}(z)k^+_{2}(z)^{-1},\\
{\cal E}(k^\pm_i(z))&=&1,\quad {\cal E}(X^\pm(z))=0.
\label{Rel5}
\ea
In the above  comultiplication  
the central elements $q^{\pm {c\over 2}}$ act nontrivially 
on the modules 
$V_1$ or $V_2$ appearing in a tensor product $V_1\otimes V_2$ 
according to whether they are written as 
$q^{\pm {c\over 2}\otimes 1}$ or
$q^{\pm 1\otimes {c\over 2}}$ respectively. Since both the above
realizations of $\ugl$ are equipped with Hopf algebra structures
it is therefore important to investigate whether they are
Hopf-algebra isomorphic to each other. As we will show later
the algebra isomorphism between them found in Ref. \cite{DiFr93}
is not a Hopf algebra isomorphism.

\subsection{\bf An algebra isomorphism}

Here we briefly recall the algebra isomorphism $\varphi$ 
between the above two realizations \cite{DiFr93}. It is based
on the uniqueness of the Gauss decomposition of the 
$L^\pm(z)$ matrix operators.
In Ref. \cite{DiFr93} it is shown that the 
 map $\varphi:\quad U(R)\ra\ugl$ defined
by 
\bac
\varphi(e^+(z_-)-e^-(z_+))&=&X^+(z),\\
\varphi(f^+(z_+)-f^-(z_-))&=&X^-(z),\\
\varphi(k^\pm_i(z))&=&k^\pm_i(z),
\label{Rel6}
\ea 
where the currents $e^\pm(z)$ and $f^\pm(z)$ being uniquely
defined by the following Gauss decompositions of  $L^\pm(z)$:  
\be 
L^\pm(z)=  
{\scriptsize\left(\begin{array}{ll}1&0\\e^{\pm}(z)&1\end{array} 
\right)} 
{\scriptsize\left(\begin{array}{ll}k_1^{\pm}&0\\0&k^{\pm}_2 
\end{array}\right)} 
{\scriptsize\left(\begin{array}{ll}1&f^{\pm}(z)\\0&1\end{array} 
\right)}= 
{\scriptsize\left(\begin{array}{ll}k^\pm_{1}(z)& 
k_1^{\pm}(z)f^{\pm}(z)\\e^{\pm}(z)k_1^\pm(z)& 
e^{\pm}(z)k_1^{\pm}(z)f^{\pm}(z)+k_2^{\pm}(z)\end{array} 
\right)}, 
\label{Rel7}
\ee 
is an algebra isomorphism.
Here the currents $e^{\pm}(z)$, $f^{\pm}(z)$, and
 $k^{\pm}_i(z)$ have all Taylor expansions in $z$ or $z^{-1}$
as: 
\bac 
e^{+}(z)&=&\sum_{n=1}^\infty e^+_n z^{-n},\quad\quad
e^{-}(z)=\sum_{n=0}^\infty e^-_{-n} z^{n}, \\ 
f^{+}(z)&=&\sum_{n=0}^\infty f^+_n z^{-n},\quad\quad
f^{-}(z)=\sum_{n=1}^\infty f^-_{-n} z^{n}, \\ 
k^\pm_i(z)&=&\sum_{n=0}^\infty k^{\pm}_{i,\mp n} z^{\pm n}. 
\label{Rel8}
\ea

Now that the algebra isomorphism between the two realizations
is settled, it is natural to try to complete this 
investigation at the
level of a Hopf algebra isomorphism because, after all, 
the major
importance of these two realizations is their being equipped
with Hopf algebra structures.

\section{A Hopf algebra isomorphism }

In this section, we find a new Hopf algebra structure
for the Drinfeld realization such that $\varphi$ becomes
 a Hopf algebra isomorphism.

{\bf Definition  3.1} A homomorphism $\psi$: 
$(A,m,\eta,\Delta,S,{\cal E})
\ra (B,m^\prime,\eta^\prime,\Delta^\prime,
S^\prime,{\cal E}^\prime)$
 is a Hopf algebra isomorphism  if $\psi$ is an algebra 
isomorphism 
between $(A,m,\eta)$ and $(B,m^\prime,\eta^\prime)$\break 
satisfying
\bac
(\psi \otimes \psi)\circ\Delta&=& 
\Delta^\prime \circ\psi,\\
\psi\circ S&=&S^\prime \circ \psi,\\
{\cal E}&=&{\cal E^\prime}\circ \psi.
\label{Rel9}
\ea

Note that $\psi$ determines uniquely the Hopf algebra 
structure for $B$
if the Hopf algebra structure of $A$ is given. Moreover, we have 
the following existence result:

{\bf Proposition 3.2} {\em Let $A$ and $B$ be two associative 
algebras with
unit together with an algebra isomorphism $\psi: 
A\longrightarrow B$. 
If there exists a Hopf algebra structure on $A$, then $B$ 
can be endowed with
a Hopf algebra structure such that $\psi$ becomes a Hopf 
algebra isomorphism.}

Note that in the proposition we do not assume that $B$ has a 
priori a Hopf
algebra structure.

If such a Hopf algebra structure for $B$ exists 
then relations 
(\ref{Rel9}) must be satisfied. Then we set
\bac
&& \Delta'=(\psi\otimes\psi)\circ\Delta\circ\psi^{-1},\\
&& S'=\psi\circ S\circ\psi^{-1},\\
&& {\cal E'}={\cal E} \circ \psi^{-1}.
\label{Rel10}
\ea
By construction, $\Delta'$ is an algebra homomorphism 
$B\longrightarrow
B\otimes B$, $S'$ is an algebra antihomomorphism 
$B\longrightarrow B^{op}$
and $\cal E'$ is an algebra homomorphism 
$B\longrightarrow  \C$.
And using the fact that $\psi$ is an algebra isomorphism and 
that 
$\Delta$, $S$ and $\cal E$ make  $A$ a Hopf algebra, 
one can check that $\Delta'$, $S'$ and $\cal E'$ 
satisfy the axioms of a Hopf
algebra for $B$.

As for the isomorphism $\varphi$ of II.3, it is clear that it 
does not
preserve the Hopf algebra structure. Indeed, let us  consider
$k^\pm_1(z)$ in both 
$\ugl$ and $U(R)$:
in $\ugl$ we have $S'(k^\pm_1(z))=k^\pm_1(z)^{-1}$, but in 
$U(R)$, using $S(L^\pm(z))=L^\pm(z)^{-1}$, we find
$S(k^\pm_1(z))=k^\pm_1(z)^{-1}+f^\pm(z)k_2^\pm(z)^{-1}e^\pm(z)$.
Since $\varphi(k_1^\pm(z))=k_1^\pm(z)$, 
 $\varphi\circ S\neq S'\circ \varphi$.

Now using the proposition, we establish a different Hopf 
algebra structure for
the Drinfeld realization so that $\varphi$ becomes a Hopf 
algebra isomorphism.
We find closed formulas for $\Delta'$, $S'$ and $\cal E'$ 
leading to a new Hopf algebra structure
for $\ugl$. This is possible because $k_1^\pm(z)$ and
$k_2^\pm(z)$ are invertible using the fact that $l_{ii}^\pm$ 
are invertible
and using the Gauss decomposition.

{\bf Theorem 3.3} The isomorphism $\varphi:\quad
 U(R(z))\ra \ugl$
 is a Hopf algebra isomorphism if the Drinfeld realization of
 $\ugl$ is equipped
with the following new Hopf algebra structure:
\bac
\Delta(k^\pm_1(z))&=&k^\pm_1(z_1^\pm)\otimes
k^\pm_1(z_2^\mp)+k^\pm_1(z_1^\pm)f^\pm(z_1^\pm)\otimes 
e^\pm(z_2^\mp)
k^\pm_1(z_2^\mp),\\
\Delta(k^\pm_2(z))&=&\sum_{n=0}^\infty (-1)^n 
k^\pm_2(z_1^\pm)f^\pm(z_1^\pm)^n\otimes e^\pm(z_2^\mp)^n
k^\pm_2(z_2^\mp),\\
\Delta(e^\pm(z))&=&e^\pm(z_1^\pm)\otimes 1+ 
\sum_{n=0}^\infty (-1)^{n} 
k^\pm_2(z_1^\pm)f^\pm(z_1^\pm)^n k^\pm_1(z_1^\pm)^{-1}\otimes
e^\pm(z_2^\mp)^{n+1},\\
\Delta(f^\pm(z))&=&1\otimes f^\pm(z_2^\mp)+ 
\sum_{n=0}^\infty (-1)^{n} f^\pm(z_1^\pm)^{n+1}\otimes
k^\pm_1(z_2^\mp)^{-1}e^\pm(z_2^\mp)^n k^\pm_2(z_2^\mp),\\
S(k^\pm_1(z))&=&k^\pm_1(z)^{-1}+
f^\pm(z) k^\pm_2(z)^{-1} e^\pm(z),\\
S(k^\pm_2(z))&=&k^\pm_2(z)^{-1}-
f^\pm(z) k^\pm_2(z)^{-1}
\bigl\{\sum_{n=0}^\infty (-1)^n (k^\pm_1(z)f^\pm(z)
k^\pm_2(z)^{-1}
 e^\pm(z))^n\bigr\}\\
&&\cdot k^\pm_1(z)k^\pm_2(z)^{-1}e^\pm(z),\\
S(e^\pm(z))&=&
-\bigl\{\sum_{n=0}^\infty (-1)^n (k^\pm_1(z)f^\pm(z)
k^\pm_2(z)^{-1}
 e^\pm(z))^n\bigr\}k^\pm_1(z)k^\pm_2(z)^{-1}e^\pm(z),\\
S(f^\pm(z))&=&
-f^\pm(z) k^\pm_2(z)^{-1}
\bigl\{\sum_{n=0}^\infty (-1)^n (f^\pm(z)
k^\pm_2(z)^{-1}
 e^\pm(z)k^\pm_1(z))^n\bigr\}k^\pm_1(z),\\
{\cal E}(k^\pm_i(z))&=&1,\quad\quad
{\cal E}(e^\pm(z))={\cal E}(f^\pm(z))=0.
\label{Rel11}
\ea  
For the purpose of simplifying 
the notation we have set 
$z_1^\pm=zq^{\pm{c\over 2}\otimes 1}$ and  
$z_2^\pm=zq^{\pm 1\otimes {c\over 2}}$. The actions of $\Delta$,
$S$ and ${\cal E}$ on the currents $X^\pm(z)$ 
of $\ugl$ are easily derived
from the above expressions and relations (\ref{Rel6}).

We have found  a method
for extending an algebra isomorphism to a  Hopf algebra 
isomorphism when one of the algebras is endowed with a Hopf
structure. 
It  would be
interesting to address   
the existence of a Hopf algebra isomorphism 
between  $U(R(z))$ and $\ugl$ with  the Drinfeld
comultiplication.
Moreover, the 
same question 
could be raised for the Drinfeld-Jimbo definition of 
$\ugl$ by means
of Chevalley generators \cite{Dri85,Jim85}. 
There also Drinfeld found an 
algebra isomorphism
between the latter algebra and 
the loop realization of $\ugl$. However,  
the coalgebra  structure of the latter realization 
is not explicit.

\section{Acknowledgments}

The work of A.H.B.  is supported by the NSF Grant \# PHY9309888.
A.H.B. would like to thank  Y. Saint-Aubin and H. Aurag for
stimulating discussions. 
\end{document}